\documentclass[12pt,reqno]{amsart}
\usepackage{amsmath,amsthm}
\usepackage{amssymb}
\usepackage{verbatim}
\newtheorem{thm}{Theorem}
\newtheorem{lemma}[thm]{Lemma}
\newtheorem{cor}[thm]{Corollary}
\newtheorem{prop}[thm]{Proposition}
\newtheorem{definition}[thm]{Definition}

\usepackage{todonotes}

\newcommand{\R}{\mathbb{R}}
\newcommand{\Rn}{\mathbb{R}^n}
\newcommand{\dif}[0]{\ensuremath{\,\mathrm{d}}}

\DeclareMathOperator*{\esssup}{ess\,sup}

\def\vint_#1{\mathchoice%
          {\mathop{\kern 0.2em\vrule width 0.6em height 0.69678ex depth -0.58065ex
                  \kern -0.8em \intop}\nolimits_{\kern -0.4em#1}}%
          {\mathop{\kern 0.1em\vrule width 0.5em height 0.69678ex depth -0.60387ex
                  \kern -0.6em \intop}\nolimits_{#1}}%
          {\mathop{\kern 0.1em\vrule width 0.5em height 0.69678ex
              depth -0.60387ex
                  \kern -0.6em \intop}\nolimits_{#1}}%
          {\mathop{\kern 0.1em\vrule width 0.5em height 0.69678ex depth -0.60387ex
                  \kern -0.6em \intop}\nolimits_{#1}}}

\renewcommand{\paragraph}[1]{\noindent {\bf #1} }

\begin{document}
\title{SHADOWS OF INFINITIES}

\author[Kuusi]{Tuomo Kuusi}
\address{Department of Mathematics and Systems Analysis, Aalto University, PO~Box~11100, 00076 Aalto}
\email{tuomo.kuusi@aalto.fi}

\author[Lindqvist]{Peter Lindqvist}
\address{Department of Mathematics and Statistics, Norwegian University of Science and Technology, 7491 Trondheim, Norway}
\email{lqvist@math.ntnu.no}

\author[Parviainen]{Mikko Parviainen}
\address{Department of Mathematics and Statistics, University of
Jyv\"askyl\"a, PO~Box~35, FI-40014 Jyv\"askyl\"a, Finland}
\email{mikko.j.parviainen@jyu.fi}

\maketitle

\begin{abstract}We study unbounded ''supersolutions'' of the Evolutionary $p$-Laplace equation with slow diffusion. They are the same functions as the viscosity supersolutions. A fascinating dichotomy prevails: either they are locally summable to the power $p-1+\tfrac{n}{p}-0$ or not summable to the power $p-2+0.$
\end{abstract}

\section{Introduction}

Our object is the \emph{unbounded} supersolutions of the Evolutionary $p$-Laplace Equation
\begin{equation}
\label{equation}
\frac{\partial v}{\partial t}\:=\;\nabla \! \cdot \! (|\nabla v|^{p-2}\nabla v), \qquad \qquad 2 < p < \infty,
\end{equation}
in a domain $\Omega_T = \Omega \times (0,T),$ where $\Omega$ is a connected open domain in $\Rn.$ Here $v = v(x_1,x_2,\cdots,x_n,t)$ and $\nabla v = (\tfrac{\partial v}{\partial x_1},\tfrac{\partial v}{\partial x_2},\cdots,\tfrac{\partial v}{\partial x_n}).$ In the literature supersolutions are usually treated as weak supersolutions to the equation, but we are interested in a much wider class of functions. Our ''supersolutions'' are defined at each point in $\Omega_T,$ are lower semicontinuous, and obey a comparison principle with respect to the solutions of the equation. (There is no assumption about $\nabla v.$) If such a supersolution, in addition, is finite in a dense subset of $\Omega_T$, we call it a \emph{$p$-superparabolic function}\footnote{They were introduced in \cite{KL} under this name, but  \emph{$p$-supercaloric functions} is more consistent.}. The $p$-superparabolic  functions, now defined as in Potential Theory are, incidentally, the same functions as the \emph{viscosity supersolutions} of the Evolutionary $p$-Laplace Equation, cf. \cite{JLM}. They appear in obstacle problems and are relevant for the Perron method, see \cite{KL}.

The three cases $1 < p < 2$ (\textsf{fast diffusion}), $p = 2$ (\textsf{the Heat Equation}) and $2 < p < \infty$ (\textsf{slow diffusion}) are very different. We shall treat only the slow diffusion case $p > 2.$ In this case disturbances propagate, as it were, with finite speed. But, as we shall see, ''infinite values'' seem to propagate with infinite speed. We have detected some fascinating phenomena, which are totally absent from the linear theory.

We are interested in the set of points where ''$v(x,t) = \infty$'', the so-called \emph{infinities}. (We do not want to call them poles.) Their definition is delicate. There are several possibilities, but, first of all, the right definition must agree with the concept in the stationary case $\nabla \!\cdot \! (|\nabla u|^{p-2}\nabla u),\; u = u(x).$ The following two sets of infinities
\begin{align*}
\Xi^{\perp} &= \{(x_0,t_0)|\,\underset{(x,t)\to(x_0,t_0+)}{\lim}v(x,t) = + \infty\}\\
\Xi^{\downarrow} & = \{(x_0,t_0)|\,\underset{t \to t_0+}{\lim}v(x_0,t) = + \infty\}
\end{align*}
are of interest for a  $p$-superparabolic function $v,$ but in principle one could consider any set $\Xi$ such that $\Xi^{\perp} \subset \Xi \subset\Xi^{\downarrow}.$  In $\Xi^{\perp}$ the limit is taken via neighbourhoods of the type 
$$|x-x_0| < \rho,\quad t_0 < t < t_0 + \delta.$$
In $\Xi^{\downarrow}$ only the time variable moves. It is of utmost importance that the limits are determined only by the future times $t > t_0,$ while the past and present times $t \leq t_0$ are totally excluded from the definitions of $\Xi^{\perp}$ and $\Xi^{\downarrow}.$ This is in striking contrast to the actual pointwise value of the function, which can always be determined by only the past:
$$v(x_0,t_0) = 
\liminf_{\substack{(x,t) \to (x_0,t_0)\\ t < t_0}} v(x,t).$$ See \cite{KL1}.
 Therefore, it may so happen that $v(x_0,t_0) < \infty,$ although $(x_0,t_0) \in \Xi^{\perp}.$ This feature is  not easily dismissed. Nonetheless, we call $(x_0,t_0)$ a \emph{point of infinity} for $v$, or, just an \emph{infinity}. (In the stationary case they are \emph{poles}.)

At this stage, we interrupt our tale by introducing the celebrated Barenblatt solution
\begin{equation}
\label{Baren}
\mathfrak{B}(x,t) = 
\begin{cases} t^{-\tfrac{n}{\lambda}}\left[C-\frac{p-2}{p}\lambda^{\tfrac{1}{1-p}}\left(\frac{|x|}{t^{1/\lambda}}\right)^{\frac{p}{p-1}}\right]_{+}^{\tfrac{p-1}{p-2}},\quad \text{when} \quad t > 0\\
0,\quad \text{when} \quad t \leq 0
\end{cases}
\end{equation}
found in 1951, cf.\ \cite{B}. Here $\lambda = n(p-2)+p.$ It is a solution of the Evolutionary $p$-Laplace Equation, except at the origin $x=0, t=0.$ It is a $p$-superparabolic function in the whole $\Rn\times \R,$ where it satisfies the equation
$$\frac{\partial \mathfrak{B}}{\partial t} - \nabla\cdot(|\nabla \mathfrak{B}|^{p-2}\nabla  \mathfrak{B}) = c \delta$$
in the sense of distributions ($\delta = $ Dirac's delta). Note carefully that due to the requirement of semicontinuity, $\mathfrak{B}(0,0) = 0$ and not $=\infty$ at the 
point $(0,0)$ of infinity. Now $(0,0) \in \Xi^{\downarrow}$ but  $(0,0) \not \in \Xi^{\perp}.$ ---In passing, we cannot resist mentioning that even for the Heat Equation
$$\frac{\partial v}{\partial t} = \Delta v$$
a similar situation appears with the fundamental solution
\begin{equation*}
\mathfrak{W}(x,t) =
\begin{cases}
\dfrac{1}{(4\pi t)^{\frac{n}{2}}}e^{-\frac{|x|^2}{4t}},\quad \text{when} \quad t > 0\\
0, \quad \text{when} \quad t \leq 0.
\end{cases}
\end{equation*}
Now $\mathfrak{W}(0,0) = 0$ while $\lim_{t \to 0+}\mathfrak{W}(0,t) = \infty.$ In classical Potential Theory, one often introduces an auxiliary supercaloric function in order to include $\{(0,0)\}$ among the ''polar sets'', see \cite{W}. Such an awkward procedure is not natural for $p > 2.$ In the non-linear theory the presence of the original $p$-superparabolic function is central.

In order to proceed, we recall that the $p$-superparabolic functions were required to be finite in a dense subset. Although, this at least excludes ''supersolutions'' that are identically infinite during some time interval, arbitrarily fast growth is still  posssible. For example, there are $p$-superparabolic functions of the form
\begin{equation*}
v(x,t) =
\begin{cases}
u(x)e^{+\frac{1}{(p-2)(t-t_0)}}, \quad \text{when}\quad t > t_0\\
0, \quad \text{when} \quad t \leq t_0
\end{cases}
\end{equation*}
where $u(x) > 0$ in $\Omega.$ Notice that here the set $\Xi^{\perp} = \Omega \times \{t_0\}.$ As we shall see, the property that the infinities occupy the whole $\Omega$  at some  time $t_0$ is a typical phenomenon for a class of $p$-superparabolic functions.

An important result is that the $p$-superparabolic functions $v:\,\Omega_T \to (-\infty,\infty]$ are of two different kinds:
\begin{enumerate}
\item[Class $\mathfrak{B}$]\qquad \quad $v \in L_{loc}^{p-1+\frac{p}{n}-\varepsilon}(\Omega_T)$ for each $\varepsilon > 0$\\
\item[Class $\mathfrak{M}$]\qquad $v \not \in L_{loc}^{p-2+\varepsilon}(\Omega_T)$ for each  $\varepsilon > 0$\\
\end{enumerate}
%
%
%
There is no third possibility\footnote{If the assumption about finite values in a dense set is abandoned, one has a further class III, see Section 6. It is outside the scope.}. The classes are not empty. The void gap $$\bigl(p-2,\,p-1+\frac{p}{n}\bigr)$$ is remarkable, to say the least. In other words, if $v \in L_{loc}^s(\Omega_T)$ for some $s > p-2,$ then $v \in L_{loc}^q(\Omega_T)$ whenever $q < p-1+\tfrac{p}{n}$ (Lemma \ref{gap}). The functions of class  $\mathfrak{B}$  have the important property\footnote{In principle, this was settled in  \cite{KL1} and \cite{KL2}, but the existence of class $\mathfrak{M}$ was, unfortunately, overlooked there.} that their gradients $\nabla v$ exist in Sobolev's sense and
$$ \nabla v \in L_{loc}^q(\Omega_T)\qquad \text{whenever} \qquad q < p-1+ \frac{1}{n+1}$$
(Theorem \ref{gradient}).
As a consequence, there exists a Radon measure $\mu \geq 0$ depending, of course, on $v$
such that
\begin{equation}
\label{measure}
\frac{\partial v}{\partial t} - \nabla \! \cdot \! (|\nabla v|^{p-2}\nabla v) = \mu
\end{equation}
in the sense of distributions. The $p$-superparabolic functions of class $\mathfrak{B}$ have a well established theory, described in  \cite{BDGO} for example. See also \cite{KLP}.

The   functions in class $\mathfrak{M}$ seem to have few good properties. First, they do not induce a Radon measure. Second, strictly speaking, their Sobolev derivative $\nabla v$ does not exist. Thus it is important to achieve simple criteria to detect functions of class $\mathfrak{M}$. Fortunately, their sets of infinities always contain a whole time slice $t=t_0,$ i.e., $v(x,t_0+) \equiv \infty$ when $x \in \Omega.$  This cannot happen for the $\mathfrak{B}$-class. The following criterion also assures that if there are too many infinities inside the domain at the same time, they have to touch the lateral  boundary. They cast their shadows on the boundary.

\begin{thm}[Theorema Infinitorum] 
A $p$-superparabolic function is of class $\mathfrak{M}$ if and only if there is a time $t_0$ such that
$$\Xi^{\downarrow}(t_0) \equiv \{(x,t_0)\in \Xi^{\downarrow}\} = \Omega \times \{t_0\}.$$
Moreover, if $\Xi^{\downarrow}(t_0)$ has positive $n$-dimensional Lebesgue measure\footnote{Recall that $\Omega \subset \Rn.$}, then $\Xi^{\downarrow}(t_0) =  \Omega \times \{t_0\}.$  The same holds for 
 $\Xi^{\perp}(t_0).$
\end{thm}

Recall that always    $\Xi^{\perp}(t_0) \subset \Xi^{\downarrow}(t_0).$  A peculiarity, which we have found, appears when $\Omega$ is the whole space $\Rn.$ In class $\mathfrak{M}$ there are no  non-negative $p$-superparabolic functions defined in $\Rn \times (0,T),$ see Theorem \ref{none}. At first sight, their absence is surprising.

Our method is based on the \emph{bounded} $p$-superparabolic functions 
$$v_k = v_k(x,t) = \min\{v(x,t),k\}.$$ 
Bounded $p$-superparabolic functions belong to the natural Sobolev space\\ $L^p_{loc}(0,T,W^{1,p}_{loc}(\Omega))$ and are weak supersolutions of the equation (\ref{equation}), cf. \cite{KL1}, \cite{LM1}. Thus \emph{a priori} estimates like the Caccioppoli inequalities are available. A convenient version of an inequality of Harnack's type, given in \cite{K2}, is needed in our work. It is valid for non-negative $p$-superparabolic functions that are bounded, in particular for the above $v_k$'s. (See Lemma \ref{Harnack} below.)

 Let us first mention a few selected results for class $\mathfrak{B}$:

\begin{thm}[Class $\mathfrak{B}$] For a $p$-superparabolic function $v:\, \Omega_T \to (- \infty,\infty]$ the following conditions are equivalent:
\begin{enumerate}
\item[(i)] \quad$v \in L_{loc}^s(\Omega_T)\qquad \text{for some}\quad s > p-2.$
\item[(ii)] \quad$\nabla v \quad \text{exists and} \quad  \nabla v \in L_{loc}^q(\Omega_T) \quad \text{whenever} \quad q < p-1+ \frac{1}{n+1}.$
\item[(iii)] \quad When $\delta>0,$
$$\underset{0<t<T-\delta}{\esssup}\int_{D}\!|v(x,t)|\,\dif x\quad < \quad \infty,\quad \text{if}\quad D \subset \subset\Omega.$$
\item[(iv)] \quad The $n$-dimensional measure $|\Xi^{\downarrow}(t)| = 0$ for each $0 < t < T$
\item[(v)] \quad It never happens that 
$$
\lim_{\substack{(y,t) \to (x,t_0)\\t>t_0}} v(y,t) =  \infty\quad \text{for all}\quad x \in \Omega.$$
\end{enumerate}
\end{thm}

\smallskip

Condition (ii) has the important implication that there exists a non-negative Radon measure $\mu$, depending on $v,$ such that the equation
\begin{equation}
\label{radon}
\int_{0}^{T}\!\!\int_{\Omega}\Bigl(-v\frac{\partial \varphi}{\partial t}+ \langle|\nabla v|^{p-2}\nabla v, \nabla \varphi \rangle \Bigr) \dif x \dif t = \int_{\Omega_T}\!\varphi\,\dif \mu
\end{equation}
holds for all test functions $\varphi \in C_{0}^{\infty}(\Omega_T).$ In other words, equation (\ref{measure}) holds in the sense of distributions, cf. \cite{KLP}.

Then we characterize the $\mathfrak{M}$ class:
 
\begin{thm}[Class $\mathfrak{M}$] For a $p$-superparabolic function $v:\, \Omega_T \to (- \infty,\infty]$ the following conditions are equivalent:
\begin{enumerate}
\item[(i)] \quad $v \not \in L_{loc}^{p-2+\varepsilon}(\Omega_T)\quad \text{for any} \quad \varepsilon > 0.$
\item[(ii)] \quad For some $\delta>0,$
 $$\underset{0<t<T-\delta}{\esssup}\int_{D}\!|v(x,t)|\,\dif x\quad = \quad \infty\quad \text{when}\quad  D \subset \subset \Omega,\quad |D|>0.$$
\item[(iii)] \quad There is $t_0$ such that the $n$-dimensional measure $|\Xi^{\downarrow}(t_0)| > 0.$
\item[(iv)] \quad    There is $t_0$ such that  $\Xi^{\downarrow}(t_0) = \Omega \times \{t_0\}.$
\item[(v)] \quad At some point $(x_0,t_0),$
$$
\liminf_{\substack{(y,t) \to (x_0,t_0)\\t>t_0}} \bigl(v(y,t)t^{\frac{1}{p-2}}\bigr) > 0.$$
\end{enumerate}
%
\end{thm}

From this one can read off a simple sufficient condition to guarantee that a $p$-superparabolic function belongs to class $\mathfrak{B}$. It is clear that a function which is bounded near the boundary cannot belong to class $\mathfrak{M}$. More precisely, if
$$\limsup_{(x,\tau)\to(\xi,t)}v(x,\tau)\: <\; \infty\quad\text{for all}\quad (\xi,t) \in \partial \Omega \times [0,T)$$
then $v$ is of class $\mathfrak{B}$.
A further result, which we find  astonishing, is that $\Xi^{\perp}$ cannot contain a portion with positive area of any other hyperplane intersecting $\Omega_T$ than those of the type $t= $Const.. 

\begin{prop}
\label{hyperplane}
If\,  $\Xi^{\perp}$ contains a portion with positive area of the hyperplane 
$$t = \langle a,x\rangle + \alpha$$
then $a = \overline{0}.$
\end{prop}
The $p$-superparabolic functions of class $\mathfrak{M}$ do not induce a $\sigma$-finite measure $\mu$ and are therefore beyond the scope of most articles devoted to the Evolutionary $p$-Laplace Equation. The fatal feature is that the possibility 
$$\iint_{K} \!|\nabla v|^{p-1} \dif x \dif t   = \infty,\qquad K \subset \subset \Omega_T,$$
cannot be avoided, in which case equation (\ref{radon}) does not make sense. Moreover, $\nabla v$ does not exist in Sobolev's sense. However, if $v\geq 1$ we know from \cite{KL2}, Theorem 4.3, that $\nabla \log v$ exists, and the Caccioppoli estimate
\begin{gather}
\label{logarithm}
\int_0^T\!\!\int_{\Omega}|\nabla \log v|^p\zeta^p\,\dif x \dif t \leq c\int_0^T\!\!\int_{\Omega}v^{2-p}\Bigl\vert\frac{\partial \zeta^p}{\partial t}\Bigr\vert \dif x \dif t \nonumber\\
+ c \int_0^T\!\!\int_{\Omega}|\nabla \zeta|^p\,\dif x  \dif t
\end{gather}
holds.

\medskip

\paragraph{Infinite initial values for solutions.} There are interesting consequences for  the solutions of the Cauchy-Dirichlet problem in $\Omega_T$:
\begin{equation}
\begin{cases}
\frac{\partial u}{\partial t}\;=\;\nabla\! \cdot\! \bigl(|\nabla u|^{p-2}\nabla u \bigr)\\
u(x,0) \;=\;g(x),\quad \text{when} \quad x \in \Omega,
\end{cases}
\end{equation}
where $0\,\leq g(x)\,\leq \infty,$ if infinite initial values are prescribed. The lateral boundary values are not essential now. Let us suppose that $u$ is a weak solution in $\Omega_T$ and that $u \geq 0$ in $\Omega_T$, see Definition 5 in Section 2. We assume that the initial values are infinite in a set $E\;\subset\;\Omega$:
\begin{equation}
\label{initi}
\lim_{t \to 0+}u(x,t)\;=\;+\infty\quad\text{for all}\quad x \in E.
\end{equation}
The following results come  from our study:
\begin{itemize}
\item If the $n$-dimensional measure of $E$ is strictly positive, then $E\,=\,\Omega.$
\item There exist solutions, if  $E\,=\,\Omega$ and $\Omega$ is bounded.
\item If $\Omega\,=\, \Rn$  and the measure of $E$ is positive, then there is no solution.
\end{itemize}
To spell it out, the requirement that
$$\limsup_{t\to 0+}u(x_0,t)\;< \; \infty$$
at some point $x_0$ is incompatible with the condition that (\ref{initi}) holds in a set $E$ of positive measure. If we replace $\Omega_T$ with a domain like
$$\Upsilon \;=\;\{(x,t)\,|x\in \Omega,\,\Psi(x) < t < T \},$$
where $\Psi\,=\,\Psi(x)$ is a smooth function, then the corresponding initial condition
$$\lim_{t\to\Psi(x)+}u(x,t)\;=\;\infty  \quad\text{at every}\quad x\in \Omega$$
is impossible, except when $\Psi(x) =$ constant. Thus we are back to the  space$\times$time cylinders. --- We hope to return to this matter in a later  work.

\section{Preliminaries}

\sloppy
We begin with some standard notation. We consider an open domain $\Omega$ in $\Rn$ and denote by $L^p(t_1,t_2;W^{1,p}(\Omega))$ the Sobolev space of functions $v = v(x,t)$ such that for almost every $t,\, t_1 \leq t \leq t_2,$ the function $x \mapsto v(x,t)$ belongs to
$W^{1,p}(\Omega)$ and
$$\int_{t_1}^{t_2}\!\!\int_{\Omega}\bigl(|v(x,t)|^p + |\nabla v(x,t)|^p\bigr)\dif x \dif t \;< \;\infty,$$
where $\nabla v = (\tfrac{\partial v}{\partial x_1},\cdots,\tfrac{\partial v}{\partial x_2}).$ The definition of the local space $L^p(t_1,t_2;W^{1,p}_{loc}(\Omega))$ is analogous. The space $L^p_{loc}(t_1,t_2;W^{1,p}_{loc}(\Omega))$ is also used.

\begin{definition}
Let $u \in L^p(t_1,t_2;W^{1,p}(\Omega)).$ Then $u$ is a weak \emph{solution} of the Evolutionary $p$-Laplace Equation in $\Omega \times (t_1,t_2),$ if 
\begin{equation}
\label{solution}
 \int_0^T\!\!\int_{\Omega}\Bigl(-u\frac{\partial \varphi}{\partial t} + \langle|\nabla u|^{p-2}\nabla u, \nabla \varphi\rangle\Bigr) \dif x \dif t\; =\; 0
\end{equation}
whenever $\varphi \in C_0^{\infty}(\Omega \times (t_1,t_2)).$ If, in addition, $u$ is continuous, then it is called a \emph{$p$-parabolic function}. Further, we say that $u$ is a weak \emph{supersolution}, if the above integral is $\geq 0$ for all $\varphi \geq 0$ in  $C_0^{\infty}(\Omega \times (t_1,t_2)).$ If the integral is non-positive instead, we say that $u$ is a weak \emph{subsolution}.
\end{definition}

By parabolic regularity theory, a weak solution is locally H\"{o}lder continuous after a possible redefinition in a set of $n\!+\!1$-dimensional Lebesgue measure zero, see \cite{T2} and \cite{Db}. Also a weak supersolution can be made semicontinuous  through such a redefinition, cf. \cite{K1}. Then it is a $p$-superparabolic function according to the Comparison Principle below.

\begin{lemma}[Comparison Principle]
Assume that  $u$ and $v$ belong to\\ $L^p\bigl(t_1,t_2;W^{1,p}(\Omega)\bigr) \bigcap C\bigl(\overline{\Omega} \times [t_1,t_2)\bigr).$ If $v$ is a weak supersolution and $u$ a weak subsolution  in $\Omega_{t_1,t_2} = \Omega \times  (t_1,t_2)$ such that 
$$ v \geq u \quad \text{on the parabolic boundary} \quad \overline{\Omega} \times \{t_1\} \bigcup
\partial  \Omega \times  (t_1,t_2),$$
then $ v \geq u$ in the whole  $\Omega_{t_1,t_2}.$
\end{lemma}

The Comparison Principle is used to define  $p$-superparabolic functions:

\begin{definition} A function $v:\,\Omega \times  (t_1,t_2)\mapsto (-\infty,\infty]$ is called a \emph{$p$-superparabolic function} if the conditions
\begin{itemize}
\item  $v$ is lower semicontinuous
\item  $v$ is finite in a dense subset
\item  $v$ satisfies the comparison principle on each cylinder $D_{t'_1,t'_2} = D \times (t'_1,t'_2)  \subset \subset \Omega_{t_1,t_2}$: if $h \in C(\overline{D_{t'_1,t'_2}})$
  is a $p$-parabolic function in  $D_{t'_1,t'_2},$ and if $h \leq v$ on the parabolic boundary of $D_{t'_1,t'_2}$,
then $h \leq v$ in the whole $D_{t'_1,t'_2}$
\end{itemize}
are valid.
\end{definition}

We recall a fundamental result from \cite{KL2}, Theorem 1.4; see \cite{LM1} for a better proof based on infimal convolutions. See also \cite{KKP}.

\begin{thm} Let $p \geq 2.$ If $v$ is a $p$-superparabolic function that is locally bounded from above in $\Omega_T,$ then the Sobolev gradient $\nabla v$ exists and $\nabla v \in L^p_{loc}(\Omega_T).$ Moreover, $v$ is a weak supersolution.
\end{thm}

In order to derive estimates from the theorem, we need \emph{bounded} functions. The truncations $v_k = \min\{v(x,t),k\}$ are $p$-superparabolic, if $v$ is, and they are bounded from above. Thus $\nabla v_k$ is at our disposal and estimates derived from the inequality
\begin{equation}
\int_0^T\!\!\int_{\Omega}\Bigl(-v_k\frac{\partial \varphi}{\partial t} + \langle|\nabla v_k|^{p-2}\nabla v_k, \nabla \varphi\rangle\Bigr)\! \dif x \dif t \;\geq \;0
\end{equation}
where $\varphi \geq 0$ and  $\varphi \in C_0^{\infty}(\Omega_T)$ are available.
The usual Caccioppoli estimates are valid.

\begin{lemma}[Caccioppoli]
\label{Caccioppoli}
Let $p > 2.$ Assume that $v \geq 1$ is a weak supersolution in $\Omega_T.$ Then the estimate
\begin{gather*}
\int_{0}^{T}\!\!\int_{\Omega}|\nabla(\zeta v^{\frac{p-1-\beta}{p}})|^p\dif x \dif t+ \frac{1}{|\beta -1|}\underset{0<t<T}{\esssup}\int_{\Omega}v^{1-\beta}\zeta^p \dif x\\
\leq C(p)(\beta^{p-1}\! +\! \beta^{1-p})\left\{\int_{0}^{T}\!\!\int_{\Omega}v^{p-1-\beta}|\nabla \zeta|^p \dif x \dif t 
+\frac{1}{|\beta -1|}\int_{0}^{T}\!\!\int_{\Omega}v^{1-\beta}\Bigl|\frac{\partial \zeta ^p}{\partial t}\Bigr| \dif x \dif t  \right\}
\end{gather*}
holds for all $\beta > 0.$\footnote{When $\beta \approx 1$ the quantity $v^{1-\beta}/|\beta - 1|$ should be replaced by
$$\Bigl|\frac{v^{1-\beta}-1}{1-\beta}\Bigr|.$$}
For $\beta = 1$ we have
\begin{gather*}
\Bigl(\frac{p}{p-2}\Bigr)^p \int_{0}^{T}\!\!\int_{\Omega}|\zeta \nabla v^{\frac{p-2}{p}}|^p\dif x + \underset{0<t<T}{\esssup}\int_{\Omega}\zeta^p\log(v) \dif x\\
\leq C(p) \int_{0}^{T}\!\!\int_{\Omega}v^{p-2}|\nabla \zeta|^p \dif x \dif t 
+\int_{0}^{T}\!\!\int_{\Omega}\log(v)\Bigl|\frac{\partial \zeta ^p}{\partial t}\Bigr| \dif x \dif t.  
\end{gather*}
Here $\zeta \geq 0$ is an arbitrary test function in $C_0^{\infty}(\Omega_T).$
\end{lemma}

\medskip

\emph{Proof:}  A formal calculation with the test function $\phi = v^{-\beta}\zeta^p$ yields the inequality. (The cases $\beta\,>\,1$ and $\beta\,<\,1$ are different.) See \cite{Db}, \cite{KL1}, \cite{K2}.
A variant of Harnack's inequality is expedient in our present work. It is valid for supersolutions.

\begin{lemma}[Harnack] Let $p>2.$ If $v > 0$ is a lower semicontinuos weak supersolution in $B(x_0,4R)\times (0,T),$ then the inequality
\begin{equation}
\label{Harnack}
\vint_{B(x_0,R)}\!v(x,t)\dif x \leq \frac{1}{2} \Bigl(\frac{c_1R^{p}}{T-t}\Bigr)^{\frac{1}{p-2}} +c_2\, \underset{Q_{2R}}{\mathrm{ess\,inf}}\{v\},
\end{equation}
 where
\begin{gather*}
Q_{2R} = B(x_0,2R) \times \bigl(t+\frac{\tau}{2}, t+\tau \bigr),\\
\tau = \min\Bigl\{T-t,\,c_1R^p\Bigl(\vint_{B(x_0,R)}\!v(x,t)\dif x \Bigr)^{2-p}\Bigr\},
\end{gather*}
is valid at a.e. time $ t,\,0<t<T.$ Here $c_1 = c_1(n,p),\,c_2 = c_2(n,p).$
\end{lemma}

This is  Theorem 1.1 in \cite{K2}. Note that the waiting time $\tau$ depends on $t.$ The estimate is valid for the so-called Lebesgue times, as explained in \cite{K2}. We only need to know that they are dense in $(0,T).$  The convenient
notation 
$$\vint_{B(x_0,R)}\!v(x,t)\dif x \; = \; \frac{\int_{B(x_0,R)}\!v(x,t)\dif x}{\int_{B(x_0,R)}\dif x}$$
is used for the average value.

\section{Examples and Comments}

We shall illustrate the theory with several examples. We begin with a simple observation.\\

\paragraph{Extension to the past.}  If $v$ is a $p$-superparabolic function in $\Omega \times (0,T)$ and if $v \geq 0$ there, then the extended function
\begin{equation}
\label{extension}
v(x,t) =
\begin{cases}
v(x,t), \qquad \text{when} \qquad 0<t<T\\
0,\qquad\quad \text{when}\qquad t\leq 0
\end{cases}
\end{equation}
is $p$-superparabolic in  $\Omega \times (-\infty,T).$ We use the same notation for the extended function.\\

\paragraph{The stationary case.}  If $v(x,t)=u(x),$ i.e., $v$ is independent of $t,$ the equation becomes the elliptic $p$-Laplace equation
\begin{equation}
\label{elliptic}
\nabla\!\cdot\!\bigl(|\nabla u(x)|^{p-2}\nabla u(x)\bigr)\;=\;0
\end{equation}
in the domain $\Omega.$ The $p$-superparabolic functions become the \emph{$p$-superharmonic functions}, defined in \cite{L}. A typical unbounded one is the fundamental solution
$$u(x)\;=\;C_{n,p}\left|x-x_0\right|^{\frac{p-n}{p-1}},$$
if $1<p<n.$ (The function is bounded if $p>n,$ and the singularity at $x=x_0$ escapes the definition, because it is not an infinity. A singularity it is.)

The infinities can be dense in the domain. We give the example
$$v(x,t)\;=\;u(x)\;=\; \sum_{j=1}^{\infty}\,\frac{C_{j,p}}{|x-q_j|^{\frac{n-p}{p-1}}}\qquad(2<p<n)$$
where $q_1,q_2,q_3,\cdots$ are the rational points and the $C_{j,p}$'s are positive convergence factors. The function is, indeed, $p$-superharmonic in $\Rn$, see \cite{LM2}. At each rational point
$$u(q) \;=\;\lim_{x \to q}u(x) \;=\;\infty.$$
This means that $v=v(x,t)$ is a $p$-superparabolic function taking the value $\infty$ along each rational line $(q,t),\;-\infty < t < \infty.$ In this case\footnote{In fact, $\lim u = \infty$ also at some irrational points. The set of infinities (the poles) is a $G_{\delta}$ set of zero $p$-capacity.}
$$\mathbb{Q}^n\times (-\infty,\infty) \;\subset \;\Xi^{\downarrow}\;=\;\Xi^{\perp}.$$
Nonetheless, $v$ is of class $\mathfrak{B}$. In particular, $\nabla v \in L_{loc}^q(\Rn\times\R)$ whenever $q < p-1+\tfrac{1}{n+1}.$ Now, as always in the stationary case, the exponent has the  better range $q < \tfrac{n(p-1)}{n-1}$ according to \cite{L}.\\

\paragraph{The Barenblatt Solution.} This function was treated in the Introduction.\\

\paragraph{A Separable Minorant.} If $\Omega$ is a bounded regular domain, there exists a $p$-superparabolic function of the form
\begin{equation}
\label{separable}
V(x,t) =
\begin{cases}
\frac{\mathfrak{U}(x)}{(t-t_0)^{\frac{1}{p-2}}},\quad \text{when}\quad t > t_0\\
0, \qquad \qquad \text{when}\quad t \leq t_0
\end{cases}
\end{equation}
where $\mathfrak{U} \in C(\overline{\Omega})\cap W^{1,p}(\Omega)$ is a weak solution to the equation 
\begin{equation}
\label{eq:elliptic-eq}
\begin{split}
\nabla \!\cdot\!\bigl(|\nabla \mathfrak{U}|^{p-2}\nabla \mathfrak{U}\bigr)\:+\;\tfrac{1}{p-2}\,\mathfrak{U}\:=\;0
\end{split}
\end{equation}
and $\mathfrak{U}>0$ in $\Omega.$ Moreover, one can take $\mathfrak{U}|_{\partial \Omega} = 0.$ The function $V$ is $p$-parabolic, when $t > t_0.$ The solution $\mathfrak{U}$ is unique\footnote{
The corresponding solution for the Porous Medium Equation is sometimes called ''the friendly giant'', see  \cite{DK} or  page 111 in \cite{V}.}. To construct $\mathfrak{U}$, we first minimize the Rayleigh quotient
$$R(w)\;=\;\frac{\int_{\Omega}|\nabla w|^p\dif x}{\Bigl(\int_{\Omega}\!|w|^2\dif x\Bigr)^{\frac{p}{2}}}$$
among all functions $w$ in $W^{1,p}_0(\Omega),\;w\not \equiv 0.$ Since $R(|w|) = R(w),$ we may assume that $w\geq 0.$ By Sobolev's and H\"{o}lder's inequalities $R(w) \geq C(n,p,|\Omega|) > 0$ for all admissible $w.$ The direct method in the Calculus of Variations yields the existence of a minimizer $u \geq 0,\: u \not \equiv 0,$ which satisfies the Euler-Lagrange equation
$$ \nabla \!\cdot\!\bigl(|\nabla u|^{p-2}\nabla u\bigr)\;+\;\lambda \|u\|_{L^2(\Omega)}^{p-2}u\;=\;0,$$
where $\lambda > 0$ is the minimum sought for. We need a normalization. Fix $u$ so that $\|u\|_{L^2(\Omega)} =1$ and note that $Cu$ satisfies the equation
$$ \nabla \!\cdot\!\bigl(|\nabla  (C u)|^{p-2}\nabla ( Cu)\bigr)\;+\;\lambda\,C^{p-2}(Cu)\;=\;0.$$
Then choose the constant $C$ so that $\lambda C^{p-2} = \tfrac{1}{p-2}.$ The so obtained $\mathfrak{U} = Cu$ is the desired solution. By elliptic regularity theory $\mathfrak{U} \in C(\overline{\Omega}) $ and $\mathfrak{U}|_{\partial \Omega} =0.$ Finally, since \,$\nabla \!\cdot\!\bigl(|\nabla \mathfrak{U}|^{p-2}\nabla \mathfrak{U}\bigr)\;\leq 0$ \, and $\mathfrak{U}\geq 0$ in $\Omega$, Harnack's inequality for supersolutions of the elliptic $p$-Laplace equation implies that $\mathfrak{U} > 0$ in $\Omega.$ See \cite{T1}. ---We could also prescribe other non-negative boundary values for $\mathfrak{U}$, but these are less needed. In only one space dimension, a formula for $\mathfrak{U}$ is easily obtained.

The constructed function $V = V(x,t)$ is a $p$-parabolic function, when $t>t_0.$ This is a useful property, since it can serve as a minorant. The functions of class $\mathfrak{M}$ have to blow up at least as fast as $(t-t_0)^{-1/(p-2)}.$

\begin{lemma} If $v \geq 0$ is a $p$-superparabolic function in $\Omega_T$ and if
\begin{equation*}
 \lim_{\substack{(y,t) \to (x,0)\\t>0}}v(y,t)\:=\;\infty
\end{equation*}
for all $x \in \Omega,$ then 
$$v(x,t)\, \geq \, \frac{\mathfrak{U}(x)}{t^{\frac{1}{p-2}}}\quad   \text{in} \quad \Omega_T.$$
In particular, 
$$ 
 \liminf_{\substack{(y,t) \to (x,0)\\t>0}}(t^{\frac{1}{p-2}}v(x,t))\; > 0\quad  \text{in} \quad \Omega.$$
\end{lemma}

\emph{Proof:} The comparison principle yields that 

$$v(x,t)\, \geq \, \frac{\mathfrak{U}(x)}{(t+\sigma)^{\frac{1}{p-2}}},\qquad \text{in}\quad \Omega_{T-\sigma},$$
where $\sigma > 0$ is arbitrarily small. Let $\sigma \to 0.$ \qed

\medskip

A superposition of a finite number of these functions is possible. Indeed,
$$ v(x,t) \;=\;\mathfrak{U}(x)\,\sum_{j=1}^{N}\Bigl[\frac{1}{t-t_j}\Bigr]_+^{\frac{1}{p-2}}$$
is $p$-superparabolic. (This construction does not work for $N=\infty$.)

The previous Lemma gives the slowest possible growth for $p$-superparabolic functions of class $\mathfrak{M}$  . But there is no upper bound. The growth can be arbitrarily fast. We just give the example
\begin{equation}
\label{exp}
v(x,t) =
\begin{cases}
\mathfrak{U}(x)\exp\Bigl(\frac{1}{(p-2)(t-t_0)}\Bigr)\quad \text{when}\quad t > t_0\\
0\quad\qquad \qquad \qquad \qquad\;\; \text{when}\quad t \leq t_0.
\end{cases}
\end{equation}
Here $\Xi^{\downarrow} = \Omega \times \{t_0\}.$ ---One can even build a tower of exponentials to increase the terrible speed of growth.\\

\paragraph{Hyperplanes in $\Xi^{\perp}.$}
As we have seen, $\Xi^{\perp}$ and  $\Xi^{\downarrow}$ can  contain portions of planes of the form $t = t_0,\,$ so-called time slices. But, surprisingly enough, no planes like 
$$ t = \langle a,x\rangle + t_0, \qquad a \not = 0,$$
will do. Indeed, the associated ''supersolution'' would be identically $\infty,$ when
$\inf\langle a,x\rangle < t-t_0 < \sup\langle a,x\rangle.$ This is outside the realm of $p$-superparabolic functions, violating the requirement of a dense subset of finite values.

\medskip

\emph{Proof of Proposition \ref{hyperplane}:} To simplify the exposition, we first treat the case with only one space variable ($n=1$). Assume that $v$ is $p$-superparabolic in $(0,2)\times(-\infty,\infty)$ and that $\Xi^{\perp}$ contains the line segment
$$t= ax,\quad a>0,\quad 0 < x < 2.$$
This will lead to the contradiction that $v=\infty$ in too large a set. To see this, fix $0 < x_0 < 1,\, t_0 = ax_0.$ Let $k >>1.$ We claim that
$$ v(x,t)\;\geq \; k\,\frac{x-x_0}{(t-t_0+\sigma)^{\frac{1}{p-2}}}, \qquad \sigma > 0,$$
in the triangular domain
$$x_0 < x < 1,\: t > ax,\; t < a\cdot 1.$$
The claim follows from the comparison principle, because the minorant is a smooth subsolution and the inequality is obviously valid on the parabolic boundary: $x=x_0,\,t \geq t_0;\; x=1,\,t=a;\; t=ax,\, x_0 \leq x \leq 1.$ Send $k$ to $\infty.$
As a result, $v = \infty$ in the whole triangular subdomain. This is a contradiction\footnote{Needless to say, the line could be replaced by a pretty arbitrary curve, and again only the time slices $t$ = Constant are acceptable to avoid a contradiction.}. This was the case $n=1.$

The proof in several dimensions is rather similar. The equation is invariant under rotations and reflexions of the $x$-coordinates. Therefore we may assume that $a_1>0, a_2>0,\cdots,a_n>0$ in the equation
$$t=a_1x_1+a_2x_2+ \cdots +a_nx_n+t_0$$
for the plane. The function
$$u(x,t) = kt^{-\frac{1}{p-2}}x_1x_2\cdots x_n \qquad \qquad (k>0)$$
is a $p$-subparabolic function when $t>0$ and $x_1x_2\cdots x_n >0.$ This is easy to verify by direct calculation, since the function is smooth. On the parabolic boundary of the polyhedral  domain
\begin{gather*}
0 < a_1x_1+a_2x_2+ \cdots +a_nx_n < t-t_0 <1,\\
x_1 > 0, x_2 > 0,\cdots,x_n > 0
\end{gather*}
we have 
$$v(x,t) \geq k(t-t_0+\sigma)^{-\frac{1}{p-2}}x_1x_2\cdots x_n,\qquad\quad \sigma >0 $$
for the given $p$-superparabolic function $v,$ which we tacitly assume to be defined here. (The boundary consists of parts of $n\!+\!1$ planes, but the plane $t=t_0+1$ is excluded.) By the comparison principle the inequality holds in the whole polyhedral domain. Letting $k \to \infty$ we see that $v=\infty$ in an open set, which means that $v$ cannot be finite in a dense subset. This contradiction concludes our proof. \qed
\\

\paragraph{Fast Growth.}  It is easy to exhibit $p$-superparabolic functions of the form
\begin{equation*}
v(x,t) =
\begin{cases}
\mathfrak{U}(x)\Psi(t), \quad t_0 < t < T\\
0,\qquad \qquad\quad t \leq t_0,
\end{cases}
\end{equation*}
where $\mathfrak{U}$ was constructed in connection with equation \eqref{eq:elliptic-eq}. One example with $T=\infty$ was formula (\ref{exp}).\\

\paragraph{Solutions that Blow Up:} The Evolutionary $p$-Laplace Equation has solutions that blow up at a certain time. The example
$$\mathfrak{D}(x,t)\;=\;\left\{A\Bigl(\frac{T}{T-t}\Bigr)^{\frac{n(p-2)}{\lambda (p-1)}}\;+\;\Bigl(\frac{p-2}{p}\Bigr)\lambda^{-\frac{1}{p-1}}\Bigl(\frac{|x|^p}{T-t}\Bigr)^{\frac{1}{p-1}}\right\}^{\frac{p-1}{p-2}},$$
with $\lambda=n(p-2)+p$, is given in Remark 7.1 on page 331 in \cite{Db}. It is is a $p$-parabolic function in $\Rn \times (0,T).$ It blows up at the terminal point $t=T.$ As it stays, it is outside the domain, but we can extend $\mathfrak{D}$ into the future, using for example the solution (\ref{separable}). Thus
\begin{equation*}
v(x,t) =
\begin{cases}
\mathfrak{D}(x,t), \quad \qquad \text{when} \qquad t < T\\
(t-T)^{-\frac{1}{p-2}}\mathfrak{U}(x),\qquad \text{when} \qquad t \geq T
\end{cases}
\end{equation*}
is a $p$-superparabolic function in $\Omega \times (0,\infty),$ where $\Omega$ comes from the definition of $\mathfrak{U}.$ In this case $\Xi^{\perp} = \Omega \times \{T\}.$

\section{Smoothing effects}

In this section we shall consider the summability (integrability) of $p$-superparabolic functions and their gradients. To be more precise, we show that if $v \in L^{p-2+\varepsilon}_{loc}(\Omega_T)$ then $v$ is actually in class $\mathfrak{B}$. We give alternative proofs to those in \cite{KL2}.  The basic tools are the Caccioppoli inequality in Lemma \ref{Caccioppoli}  and Sobolev's inequality, written for convenience in the form
\begin{equation}
\label{sobolev}
\int_0^T\!\!\int_{\Omega}\!|\zeta w|^q\dif x \dif t\;\leq\;S^q \int_0^T\!\!\int_{\Omega}\!|\nabla (\zeta w)|^p\dif x \dif t\, \biggl\{
 \esssup_{0<t<T}\int_{\Omega}\!|\zeta w|^m\dif x\biggr\}^{\frac{p}{n}},
\end{equation}
valid for $m>0$ and $q = p+\tfrac{pm}{n}.$ Here $\zeta \in C_0^{\infty}(\Omega_T)$ is a suitable test function and 
$$ \zeta w \in L^{\infty}(0,T;L^m(\Omega)) \cap L^p(0,T;W^{1,p}(\Omega)).$$
See \cite{Db}, Proposition 3.1, Chapter 1, page 7. Since our results are local we may as well assume that the $p$-superparabolic function $v$ is $\geq 1$ in the whole $\Omega_T.$

The estimates will be obtained by iteration. At each step  of the iteration a new test function $\zeta$ has to be chosen. Typically, the domain shrinks during the procedure. Fortunately, we need only a finite number of steps. Therefore we do not keep track of the $\zeta$'s. We begin with an alternative proof of a theorem from \cite{KL2}.

\begin{thm}
\label{gap}
Let $v$ be a $p$-superparabolic function in $\Omega_T.$ If $v \in L^{p-2+\varepsilon}_{loc}(\Omega_T)$ for some $\varepsilon > 0,$ then $v \in  L^{p-1+\frac{p}{n}- \sigma}_{loc}(\Omega_T)$ for each  $\sigma>0.$
\end{thm}

\emph{Proof:}   Since $v$ is superparabolic, it is bounded  from below and thus by adding a constant, we may assume that $v\ge 1$. Fix the desired small $\sigma > 0.$ Anticipating the procedure, we try to find an index $j$ so that
$$\varepsilon \Bigl(1+ \frac{p}{n} \Bigr)^j\;=\; 1 - \frac{\sigma}{1+\frac{p}{n}}.$$
Since the assumption holds for any $\varepsilon$ smaller than the given one, we can always accomplish this. Let $\xi \in C_0^{\infty}(\Omega_T),\; 0 \leq \zeta \leq 1,$ and set $\zeta = 1$ in any chosen compact subdomain of $\Omega_T.$ The Caccioppoli estimate
\begin{align*}
&\int_0^T\!\!\int_{\Omega}\!|\nabla (\zeta v^{\frac{p-2+\varepsilon}{p}})|^p\dif x \dif t\,
+\frac{1}{\varepsilon} 
 \esssup_{0<t<T}\int_{\Omega}\!v(x,t)^{\varepsilon}\zeta (x,t)^p\dif x\\
&\leq\; \frac{C(p)}{(1-\varepsilon)^p}\biggl\{\int_0^T\!\!\int_{\Omega}\!v^{p-2+\varepsilon}|\nabla \zeta|^p\dif x \dif t\;\,+\,\frac{1}{\varepsilon}
\int_0^T\!\!\int_{\Omega}\!v^{\varepsilon}\Bigl|\frac{\partial}{\partial t}\zeta^p \Bigr|\dif x \dif t\biggr\}
\end{align*}
in Lemma \ref{Caccioppoli} with $\varepsilon = 1-\beta$
is valid when $0 < \varepsilon < 1.$ By our assumption, the right-hand side is finite. Combining this with the Sobolev inequality  (put $w =  v^{\frac{p-2+\varepsilon}{p}}$, $m=p\varepsilon/(p-2+\varepsilon)$) 
\begin{align*}
\int_0^T & \int_{\Omega}  \zeta^{p\gamma} v^{p-2+\varepsilon \gamma} \dif x \dif t \\
& \leq\; S^{p\gamma}\int_0^T \int_{\Omega} |\nabla (\zeta v^{\frac{p-2+\varepsilon}{p}})|^p\dif x \dif t
\left\{
 \esssup_{0<t<T}\int_{\Omega}\!v(x,t)^{\varepsilon}\zeta (x,t)^p\dif x\right\}^{\frac{p}{n}}
\end{align*}
where $\gamma = 1 + \tfrac{p}{n} \,>\,1,$ we obtain
\begin{equation*}
\int_0^T\!\!\int_{\Omega}\!\zeta^{p\gamma} v^{p-2+\varepsilon \gamma}\dif x \dif t\;<\;\infty.
\end{equation*}
Thus
$$v \in L_{loc}^{p-2+\varepsilon(1+\frac{p}{n})}(\Omega_T).$$

We repeat the procedure, now with $\varepsilon(1+\frac{p}{n})$ in the place of $\varepsilon$, and obtain the better exponent
$$p-2+\varepsilon(1+\frac{p}{n})\gamma \,=\,p-2+\varepsilon(1+\frac{p}{n})^2.$$
Iterating till we reach the exponent $p-2+\varepsilon(1+\tfrac{p}{n})^j,$ we can perform one final iteration, obtaining the desired exponent 
$$p-2+\varepsilon(1+\frac{p}{n})^j\gamma \,=\,p-2+(1-\frac{\sigma}{\gamma})\gamma\,=\, p-1+\frac{p}{n}-\sigma.$$
This concludes our proof, but we remark that an explicit estimate can be worked out, which we omit, since only a finite number of iterations was needed.\qed

\medskip

In the next theorem from \cite{KL2} it is decisive that  one can deduce that $\nabla v \in L_{loc}^{p-1}(\Omega_T),$ because this is sufficient to induce a Radon measure. For the benefit of the reader, we give a proof.

\begin{thm}
\label{gradient}
Let $v \in  L_{loc}^{p-2+\varepsilon}(\Omega_T)$ be a $p$-superparabolic function in $\Omega_T.$ Then the Sobolev gradient $\nabla v$ exists and $\nabla v \in L_{loc}^{q}(\Omega_T)$ whenever $q< p-1+\tfrac{1}{n+1}.$
\end{thm}

\emph{Proof:}  The proof is the same as in \cite{L}. By \cite{KL1}, \cite{LM1} or \cite{KKP} the gradient exists.
Let $0<t_1<t_2<T$ and $K\subset \subset \Omega.$ Take $0<\beta<1.$ By the H\"{o}lder inequality
\begin{align*}
& \int_{t_1}^{t_2}\!\!\int_K\!|\nabla v|^q\dif x \dif t\\=\;& \int_{t_1}^{t_2}\!\!\int_K\!(v^{-\frac{1+\beta}{p}} |\nabla v|)^qv^{\frac{1+\beta}{p}q}\dif x \dif t\\
\leq\;& \biggl\{\int_{t_1}^{t_2}\!\!\int_K\!v^{-1-\beta}|\nabla v|^p\dif x \dif t \biggr\}^{\frac{q}{p}}\biggl\{\int_{t_1}^{t_2}\!\!\int_K\!v^{\frac{1+\beta}{p-q}q}\dif x \dif t\biggr\}^{1-\frac{q}{p}}\\
=\;&\Bigl(\frac{p}{p-1-\beta}\Bigr)^q\biggl\{\int_{t_1}^{t_2}\!\!\int_K\!|\nabla (v^{\frac{p-1-\beta}{p}})|^p\dif x \dif t\biggr\}^{\frac{q}{p}}\biggl\{\int_{t_1}^{t_2}\!\!\int_K\!v^{\frac{1+\beta}{p-q}q}\dif x \dif t \biggr\}^{1-\frac{q}{p}}.
\end{align*}
The last integral is finite  if
$$(1+\beta)\frac{q}{p-q}\:<\;p-1+\frac{p}{n}$$
by the previous theorem, and the first one by the Caccioppoli estimate, whenever $0 < \beta < 1.$ We see that any exponent 
$$q < p-1+\frac{1}{n+1}$$
is possible to reach. \qed

\medskip

\textbf{Remark:} Also the opposite implication holds: if $\nabla v \in L_{loc}^q(\Omega_T)$ when $q < p-1+\tfrac{1}{n+1},$ then $v \in L_{loc}^{p-2+\varepsilon}(\Omega_T),$ when $1 + \frac{p}{n} >\varepsilon > 0.$

\medskip

For the Barenblatt solution (\ref{Baren}) the integrals 
$$
 \esssup_{-\infty<t<\infty}\int_{\Rn}\!\mathfrak{B}(x,t)^{\alpha}\dif x$$
converge when $0 <\alpha \leq 1$ but not when $\alpha > 1.$ We have the following general result, characterizing Class $\mathfrak{B}$.

\begin{thm}
\label{alpha}
Suppose that $v$ is $p$-superparabolic in $\Omega_T.$ If $v \geq 1$ and
$$
 \esssup_{0<t<T}\int_{\Omega}\!v(x,t)^{\alpha}\dif x \: <\; \infty$$
for some exponent $\alpha >0,$ then $v \in L^s_{loc}(\Omega_T)$ whenever $s<p-1+\tfrac{p}{n}.$
\end{thm}

\medskip

\emph{Remark:} As a consequence, $$
 \esssup_{0<t<T}\int_{\Omega}\!v(x,t)^{\alpha}\dif x \: =\; \infty$$
for all exponents $\alpha >0,$ if $v$ belongs to $\mathfrak{M}.$

\medskip

\emph{Proof:}  We shall show that $v \in L^s_{loc}(\Omega_T)$ for some $s \geq p-2.$ If $s > p-2$ we are done, because Theorem \ref{gap} now applies. To be on the safe side, we first treat the case $s=p-2.$ Then $v \in L^{p-2}_{loc}(\Omega_T)$ and the Caccioppoli inequality 
\begin{align*}
\int_{0}^{T}\!\!\int_{\Omega}|\nabla& (\zeta v^{\frac{p-2}{p}})|^p\dif x \dif t \\
\leq C_1(p)&\int_{0}^{T}\!\!\int_{\Omega}\Bigl(v^{p-2}|\nabla \zeta|^p + 
\log(v)\Bigl|\frac{\partial \zeta ^p}{\partial t}\Bigr|\bigr) \dif x \dif t\\
\leq   C_2(p)&\int_{0}^{T}\!\!\int_{\Omega}v^{p-2}\Bigl(|\nabla \zeta|^p + 
\Bigl|\frac{\partial \zeta ^p}{\partial t}\Bigr|\bigr) \dif x \dif t < \infty
\end{align*}
in Lemma \ref{Caccioppoli} 
is available. In the Sobolev inequality (\ref{sobolev}) we take $w = v^{\frac{p-2}{p}}$ and $m = \tfrac{\alpha p}{p-2},$ so that $w^m = v^{\alpha}$ in the single integral. Then 
$$q= p\Bigl(1+\frac{p\alpha}{n(p-2)}\Bigr), \qquad w^q = v^{p-2+\frac{\alpha p}{p-2}}.$$
It follows that $v \in L_{loc}^{p-2+\frac{p \alpha}{n}}(\Omega_T)$ and now the summability exponent is in the range for which Theorem \ref{gap} is applicable. ---That much about the case $s=p-2.$

Next we describe an iteration, starting the procedure from some small $\alpha$ in the range $0 < \alpha< p-2.$ The Caccioppoli inequality
\begin{align*}
\int_{0}^{T}\!\!\int_{\Omega}&|\nabla (\zeta v^{\frac{\alpha}{p}})|^p\dif x \dif t \\
\leq\, &\;C_3(p)\frac{(p-1-\alpha)^{p-1}}{p-2-\alpha}\int_{0}^{T}\!\!\int_{\Omega}v^{\alpha -(p-2)}\Bigl|\frac{\partial \zeta ^p}{\partial t}\Bigr|\dif x \dif t 
\\ & \qquad + C_3(p) \int_{0}^{T}\!\!\int_{\Omega} v^{\alpha}|\nabla \zeta|^p \dif x \dif t  \\
\leq \,&\;\frac{C_4(p)}{p-2-\alpha}\int_{0}^{T}\!\!\int_{\Omega}v^{\alpha}\Bigl(|\nabla \zeta|^p + \Bigl|\frac{\partial \zeta ^p}{\partial t}\Bigr|\Bigr) \dif x \dif t \; < \; \infty 
\end{align*}
is at our disposal. (We used $v^{\alpha-(p-2)} \leq v^{\alpha}$ in the last step.)
With $w=v^{\frac{\alpha}{p}}$ and $m=p$ we can use the Sobolev inequality. Now
$$q= p\Bigl(1+\frac{p}{n}\Bigr), \qquad w^q = v^{\alpha(1+ \frac{p}{n})}.$$
It follows that $v \in L_{loc}^{\alpha(1+\frac{p}{n})}(\Omega_T).$ If $\alpha(1+\tfrac{p}{n}) < p-2$ we repeat the procedure, this time with 
$$w = v^{\alpha\frac{(1+\frac{p}{n})}{p}},\quad m = \frac{p}{1+\frac{p}{n}},\quad q = p\Bigl[1 + \frac{p}{n(1+\frac{p}{n})}\Bigr]$$
so that 
$$ w^q = v^{\alpha(1+\frac{2p}{n})}.$$
Thus the result is that $v \in L_{loc}^{\alpha(1+\frac{2p}{n})}\!(\Omega_T).$ We continue till we, sooner or later, reach an index $j$ for which
$$\alpha\Bigl(1+\frac{jp}{n}\Bigr)\;< p-2 \qquad \text{and} \qquad \alpha\Bigl(1+\frac{(j+ 1)p}{n}\Bigr)\; \geq p-2.$$
We can do one final iteration using $w = v^{\alpha(1+\frac{jp}{n})/p},\;m= \tfrac{p}{1+\tfrac{jp}{n}}.$ It follows that $v \in L_{loc}^s(\Omega_T)$ for some  $s > p-2.$
 This case was dealt with above. \qed

\section{Class $\mathfrak{M}$}

A typical $p$-superparabolic function which is not of Class $\mathfrak{B}$ is the previously constructed
$$V(x,t)\,=\,\left[\frac{1}{t-t_0}\right]_+^{\frac{1}{p-2}}\mathfrak{U}(x,t)$$ in $\Omega \times (-\infty,\infty),$ where $\Omega$ has to be  bounded. This function is not locally summable to any power $\geq p-2,$ nor is its gradient. Set, for a given function $v,$ defined in $\Omega_T$
\begin{align*}
\Xi^{\perp}(t_0)\;=\;&\bigl\{x\in\Omega|\; 
 \lim_{\substack{(y,t) \to (x,t_0)\\t>t_0}}v(y,t)\;=\;+\infty\bigr\}\\
\Xi^{\downarrow}(t_0)\;=\;&\bigl\{x\in\Omega|\; 
 \lim_{t \to t_0+}v(x,t)\;=\;+\infty\bigr\}
\end{align*}
so that
$$\Xi^{\perp}\;=\;\bigcup_{0\leq t<T}\Xi(t),\qquad \Xi^{\downarrow}\:=\; \bigcup_{0\leq t<T}\Xi^{\downarrow}(t).$$
Of course, $\Xi^{\perp}(t_0) \subset \Xi^{\downarrow}(t_0),$ but \emph{they do not have to be the same sets}, as the Barenblatt solution shows. The striking phenomenon is that if the $n$-dimensional measure $|\Xi^{\downarrow}(t_0)|\,>\,0,$ then also  $|\Xi^{\perp}(t_0)|\,>\,0.$
Before dealing with this, we need to give the following lemma about large average values. 

\begin{lemma}
\label{times}
Suppose that $v$ is a non-negative $p$-superparabolic function in $\Omega_T.$ Suppose that $B(x_0,4R) \subset  \Omega.$ If there is a sequence of ''Lebesgue times'' $t_j \to t_0, \;0<t_j<T$ such that
$$\lim_{j \to \infty} \int _{B(x_0,R)}\!v(x,t_j)\dif x \quad=\quad \infty,$$
then
$$ v(x,t)\;\geq\; \gamma\, \frac{R^{\frac{p}{p-2}}}{(t-t_0)^{\frac{1}{p-2}}},$$
when $x \in B(x_0,2R)$ and $t_0 < t < T.$ The constant $\gamma  > 0$ depends only on $n$ and $p.$
\end{lemma}

\emph{Remark:} If $t_0>0$ we do not forbid that $t_j < t_0.$

\medskip

\emph{Proof:} We aim at using Harnack's inequality (\ref{Harnack}) for the bounded supersolutions
$v_k,$ where $k$ does not have to be an integer. Now, for a fixed index $j$, by continuity the integral
$$J^k(t_j)\quad=\quad \vint_{B(x_0,R)}\!\min\{v(x,t_j),k\} \dif x$$
attains all values in the interval $[0,\vint_{}v(x,t_j)\dif x)$ when $k$ increases from $0$ to $\infty$. Let $S>t_0$ be a given number so that $S-t_0$ is small enough. Then, for $j$ large enough,
$$ \vint_{B(x_0,R)}\!v(x,t_j) \dif x\quad>\quad \left(\frac{c_1R^p}{S-t_0}\right)^{\frac{1}{p-2}}$$
 by our assumption. Determine $k_j$ (not necessarily an integer) so that
$$J^{k_j}(t_j)\;=\;\left(\frac{c_1R^p}{S-t_0}\right)^{\frac{1}{p-2}}.$$
By Harnack's inequality (\ref{Harnack}) evaluated at the Lebesgue time $t_j$ we have
\begin{gather*}
\left(\frac{c_1R^p}{S-t_0}\right)^{\frac{1}{p-2}}\;\leq \; \frac{1}{2}\left(\frac{c_1R^p}{S-t_j}\right)^{\frac{1}{p-2}}\;+\; c_2\,\inf_{Q^j_{2R}}v,
\intertext{where}
Q^j_{2R}\;=\; B(x_0,2R) \times \bigl(t_j+\frac{\tau_j}{2},t_j+\tau_j\bigr),\\
\tau_j \;=\;\min\bigl\{S-t_j,c_1R^pJ^{k_j}(t_j)^{2-p}\bigr\}\;=\;\min\{S-t_j,S-t_0\}.
\end{gather*}
Taking the limit as $j \to \infty,$ we arrive at
$$\frac{1}{2}\left(\frac{c_1R^p}{S-t_0}\right)^{\frac{1}{p-2}}\;\leq\;c_2\,\inf_{Q_{2R}} v,$$
where 
$$Q_{2R}\;=\; B(x_0,2R) \times \Bigl(\frac{t_0+S}{2},S\Bigr).$$
Then we have the inequality
$$c_2\, v(x,t)\;\geq\;\frac{1}{2}\left(\frac{c_1R^p}{S-t_0}\right)^{\frac{1}{p-2}}\;\geq\;
\frac{1}{2}\left(\frac{c_1R^p}{2(t-t_0)}\right)^{\frac{1}{p-2}},$$
when $S > t >\tfrac{S+t_0}{2}.$ By adjusting $S$ we can reach all $t$ in $(t_0,T).$ \qed

\begin{cor}
\label{onepoint}
If, at some point $(x_0,t_0),$
$$\liminf_{\substack{(x,t)\to (x_0,t_0)\\t>t_0}}\bigl((t-t_0)^{\frac{1}{p-2}}v(x,t)\bigr)\;>\;0,$$
then \,$\Xi^{\perp}(t_0)\,=\,\Omega \times \{t_0\}.$
\end{cor}

\emph{Proof:} In some small neighbourhood $|x-x_0| \leq \rho,\,t_0 < t <t_0 + \rho^p,$
we have
$$(t-t_0)^{\frac{1}{p-2}}v(x,t)\;\geq\;\varepsilon_0\;>\;0.$$
Then
$$\vint_{B(x_0,\rho)}v(x,t)\dif x\; \geq\; \varepsilon_0(t-t_0)^{-\frac{1}{p-2}}$$
as $t \to t_0+.$ The assumption in Lemma \ref{times} is fulfilled. The inclusion  
$B(x_0,2\rho)\;\subset\;\Xi^{\perp}(t_0)$ follows. We can apply  Lemma \ref{times} on any ball $B(y_0,R)$ with $B(y_0,4R) \subset \subset \Omega$ intersecting $B(x_0,2\rho)$ and conclude that also $B(y_0,R)\; \subset \;\Xi^{\perp}(t_0).$ Using a suitable chain of balls, we can see that the corollary holds.  \qed

\medskip

If there are too much infinities inside the domain, they have to touch the lateral boundary: the infinities  ''cast a shadow''. That is in the next theorem.

\begin{thm}
\label{pos}
If for some $t_0$ the $n$-dimensional measure   $|\Xi^{\downarrow}(t_0)|\,>\,0,$ then also
 $|\Xi^{\perp}(t_0)|\,>\,0.$ Actually,
$$\Xi^{\perp}(t_0)\;=\;\Xi^{\downarrow}(t_0)\;=\; \Omega \times \{t_0\}.$$
\end{thm}

\emph{Proof:} We take $t_0=0$ and select some ball $B(x_0,8R)$ in $\Omega$ so that
$$|\Xi^{\downarrow}(0)\cap B(x_0,R)|\,>\,0,$$ which is possible by the assumption. Let $k > 0.$ To each $x \in \Xi^{\downarrow}(t_0)$ there is a time $t^k_x>0$ such that
$$v(x,t) > k,\quad \text{when} \quad 0 < t < t^k_x.$$
We remark that the times $ t^k_x$ are  decreasing as $k \to \infty.$ Define the line segments
$$L^k_x\;=\; \{x\}\times (0,t^k_x)$$
and consider the projected sets
$$E^k_t\;=\; \overline{B(x_0,R)}\cap\bigl\{x|\;(x,t)\in \cup_y L_y^k\bigr\},$$
which consists of all endpoints $x \in \overline{B(x_0,R)}$ with the corresponding  $ t^k_x > t.$ The set $E^k_t$ shrinks with increasing $k.$

\textsf{Claim:}  There is a $T^k > 0$ such that 
$$|E^k_t|\;\geq\; \frac{1}{2}\,|\Xi^{\downarrow}\cap B(x_0,R)| \quad \text{when}\quad 0 < t < T^k.$$
Moreover, $T^k$ decreases when $k \to \infty.$

 Indeed,
$$\Xi^{\downarrow}(0) \cap B(x_0,R)\:=\; 
 \bigcap_{k=1}^{\infty}\bigcup_{j=1}^{\infty}E_{\frac{1}{j}}^k$$
so that
$$|\Xi^{\downarrow}(0) \cap B(x_0,R)|\;\leq\;\left|\bigcup_{j=1}^{\infty}E_{\frac{1}{j}}^k\right|\;=\; \lim_{j \to \infty}|E_{\frac{1}{j}}^k|,$$
since the sets are nested. To each $k$ there is a $j=j_k$ such that 
$$ |E_{\frac{1}{j_k}}^k|\:>\; \frac{1}{2}\,|\Xi^{\downarrow}(0) \cap B(x_0,R)|.$$
The claim follows, because $$E_{\frac{1}{j_k}}^k \subset E^k_t,\quad \text{when}\quad 0 < t < \tfrac{1}{j_k} = T^k.$$ The so defined $T^k$ are decreasing, if we select the $j_k$ to be increasing.

We select a \emph{compact} subset $\tilde{E}^t_k \subset E^t_k$ so that $|\tilde{E}^t_k|\,\geq \, \tfrac{1}{2}|E^t_k|.$ Thus
$$|\tilde{E}^t_k| \,\geq\, \frac{1}{4}\,|\Xi^{\downarrow}(0) \cap B(x_0,R)|\quad \text{when}\quad 0 < t< T^k.$$
Let $0 < t< T^k$ and $x \in \tilde{E}^t_k .$ By the semicontinuity there is a radius $r^k_{x,t} < \frac{R}{10}$ such that
$$\inf_{Q^k_{x,t}}v\;\geq \; k,\qquad Q^k_{x,t}\;=\;B(x,r^k_{x,t})\times [t,t+(r^k_{x,t})^p).$$
Obviously,
$$\tilde{E}^t_k\;\subset \bigcup_{x\in\tilde{E}^t_k}B(x,r^k_{x,t}).$$
By compactness of $\tilde{E}^t_k$ and by a simple version of Vitali's covering theorem (\cite{S}, Chapter I, paragraph 1.6) there is finite $J_k$ and \emph{disjoint} balls $B(x_j,r_j^k)$ so that 
$$\tilde{E}^t_k\;\subset \bigcup_{j=1}^{J_k}B(x_j,5r^k_j).$$
If $t< \tau < t + [\min\{r_1^k,r^k_2,\cdots r^k_{J_k}\}]^p$ then
\begin{align*}
&\int_{B(x_0,2R)}\!\min\{v(x,\tau),k\}\dif x\;\geq\; \int_{\bigcup B(x_j,r^k_j)}\min\{v(x,\tau),k\}\dif x\\
&=\;\sum_{j=1}^{J_k}\int_{B(x_j,r^k_j)}\min\{v(x,\tau),k\}\dif x\;\geq\; k\,\sum_{j=1}^{J_k}|B(x_j,r^k_j)|\\
&=\;k 5^{-n}\sum_{j=1}^{J_k}|B(x_j,5r^k_j)|\;\geq\; k5^{-n}|\tilde{E}^t_k|\;\geq\; k4^{-1}5^{-n}|\Xi^{\downarrow}(0)\cap B(x_0,R)|.
\end{align*}

 Thus
$$\int_{B(x_0,2R)}\!\min\{v(x,t),k\}\dif x\;\geq\;\frac{k}{4\cdot5^n}\,|\Xi^{\downarrow}(0) \cap  B(x_0,R)|$$
for $0 < t < T^k.$ From Lemma \ref{times} it follows that
$$
 \lim_{\substack{(y,t)\to (x,0)\\t>0}}v(y,t)\;=\;\infty\quad \text{for all}\quad x \in B(x_0,4R).$$
In other words, 
$$B(x_0,4R)\;\subset\;\Xi^{\perp}(0)\;\subset\;\Xi^{\downarrow}(0).$$
We read off from the proof that the infinities in $B(x_0,R)$ cause that the whole larger ball $B(x_0,4R)$ consists of infinities at time $t=0.$ Repeating the argument with suitable chains of balls, we can conclude that the whole $\Omega \times \{0\}$ consists of infinities. (We have assumed that $\Omega$ is connected.)\quad  \qed

\medskip

We saw that $|\Xi^{\perp}(t)| \,=\,0$ and $|\Xi^{\downarrow}(t)|\,=\,0$ simultaneously. Positive measure led to the situation with the violent behaviour described in Section 5. Yet, to complete the picture, we  need to show that, if $|\Xi^{\perp}(t)| \,=\,0$ for each $0<t<T,$ then the function belongs to Class $\mathfrak{B}$. By Theorem \ref{alpha} it is enough to establish the following.

\begin{lemma}
\label{completed}
If $|\Xi^{\perp}(t)| \,=\,0$ for each $0<t<T,$ then $v \in L^{\infty}_{loc}(0,T;L_{loc}^1(\Omega)).$
\end{lemma}

\emph{Proof:}  The \emph{antithesis} is that 
$$
 \esssup_{\varepsilon<t<T-\varepsilon}\int_{B(x_0,R)}\!v(x,t)\dif x\:=\; \infty$$
for some $R$ and $\varepsilon.$ We can extract a convergent sequence of Lebesgue times, say $t_j \to t_0,$ such that
$$
 \lim_{j \to \infty}\int_{B(x_0,R)}\!v(x,t_j)\dif x\:=\; \infty.$$
Lemma \ref{times} implies that 
$$\lim_{\substack{(y,t)\to(x,t_0)\\t > t_0}}v(y,t) \;=\;\infty\qquad \text{for all}\qquad x \in B(x_0,2R).$$
Thus $B(x_0,2R)\,\subset\, \Xi^{\perp}(t_0)$ and so $|\Xi^{\perp}(t_0)|\,>0.$ This contradiction shows that the antithesis is false. The lemma follows.  \qed

\medskip

If $\Omega\times (0,T)=\Rn\times (0,T)$, then $\Xi^\perp(0)$ is of measure zero.  
\begin{thm}
\label{none}
If  $v:\Rn \times (0,T) \to [0,\infty]$ is $p$-superparabolic, then the $n$-dimensional measure $|\Xi^{\perp}(0)|\,=\,0.$ 
\end{thm}

\emph{Proof:} Assume that  $|\Xi^{\perp}(0)|\,>\,0.$  We can regard $v$ as zero, when $t \leq 0.$ There must be a point where Corollary \ref{onepoint} applies, thus $\Xi^{\perp}(0)\,=\,\Rn \times \{0\}.$ Choose an arbitrarily large ball $B(0,R)$  and let 
$$\mathfrak{V}(x,t)\;=\; t^{-\frac{1}{p-2}}\mathfrak{U}(x)$$
be the $p$-parabolic function constructed in the unit ball $|x| < 1,$ as in formula (\ref{separable}). By scaling and comparison
$$v(x,t)\;\geq\;\Bigl(\frac{R^p}{t}\Bigr)^{\frac{1}{p-2}}\mathfrak{U}\bigl(\frac{x}{R}\bigr),\qquad t >0,\; |x| < R.$$
Let $\nu \,=\,\min_{|y|\leq 1/2}\mathfrak{U}(y)\; >\; 0.$ Then
$$v(x,t)\;\geq\; \nu \Bigl(\frac{R^p}{t}\Bigr)^{\frac{1}{p-2}} \quad \text{when} \quad |x| \leq \frac{R}{2}.$$
Letting $R \to \infty,$ we must have $v(x,t) \equiv \infty.$ The function is not finite in a dense subset. \qed

\section{The ''Outsiders''}

The $p$-superparabolic functions do not form a closed class under monotone convergence. In the stationary case, the limit of an increasing sequence of $p$-superharmonic functions is either identically infinite or a $p$-superharmonic function. For the Evolutionary $p$-Laplace Equation, the situation is not quite that simple. The limit of an increasing sequence of $p$-superparabolic functions can be a function that is identically infinite in some time intervals:
$$v(x,t)\, \equiv \infty\quad \text{when} \quad x \in \Omega,\; t_1 < t < t_2.$$
This follows from our previous considerations, because the truncations
$$\min\{v(x,t),k\}$$
are bounded $p$-superparabolic functions. It is also possible to construct examples such that
$$\Xi^{\perp}\; \supseteq \;\Omega \times \bigcup_{j} [a_j,b_j],$$
where the union of disjoint time intervals is countable.
 We can use estimate (5) to conclude that
\begin{align*}
\iint_{\Omega_T\cap\{v<\infty\}}\!|\nabla \log v|^p\zeta^p \dif x \dif t\; \leq & \; c\int_0^T\!\!\int_{\Omega}v^{2-p}\Bigl|\frac{\partial \zeta^p}{\partial t}\Bigr|\dif x \dif t\:
\\ & \; + c \int_0^T\!\!\int_{\Omega}|\nabla \zeta|^p \dif x \dif t
\end{align*}
for strictly positive $v.$ (If $v \equiv \infty$, there is nothing to say.)\\

\thanks{\emph{Acknowledgements.} This work was initiated in Parma during the workshop ''New Trends in Nonlinear Parabolic Equations'' in November 2012. The research was done at the Mittag-Leffler Institute in the autumn of 2013 under the program ''Evolutionary problems''. The final manuscript was written up at the Department of Mathematics of the University of Pittsburgh. We are pleased to thank these institutions.}

\end{document}